\documentclass[reqno,centertags]{amsart}
\usepackage{amsmath}
\usepackage{amssymb}
\usepackage{amsfonts}

\newtheorem{theorem}{Theorem}
\theoremstyle{plain}

\numberwithin{equation}{section}

\input{tcilatex}

\begin{document}
\title[]{Surfaces of coordinate finite type in the Lorentz-Minkowski 3-space}
\author{Hassan Al-Zoubi}
\address{Department of Mathematics, Al-Zaytoonah University of Jordan, P.O. Box 130, Amman, Jordan 11733}
\email{dr.hassanz@zuj.edu.jo}
\author{Alev Kelleci}
\address{Department of Mathematics, Firat University, Turkey 23100}
\email{alevkelleci@hotmail.com}
\author{Tareq Hamadneh}
\address{Department of Mathematics, Al-Zaytoonah University of Jordan, P.O.
Box 130, Amman, Jordan 11733}
\email{t.hamadneh@zuj.edu.jo}

\date{}

\begin{abstract}
In this article, we study the class of surfaces of revolution in the 3-dimensional Lorentz-Minkowski space with nonvanishing Gauss curvature whose position vector $\boldsymbol{x}$ satisfies the condition $\Delta^{III}\boldsymbol{x}=A\boldsymbol{x}$, where $A$ is a square matrix of order 3 and $\Delta^{III}$ denotes the Laplace operator of the second fundamental form $III$ of the surface. We show that such surfaces are either minimal or pseudospheres of a real or imaginary radius.
\end{abstract}

\maketitle

\section{Introduction}
Let $M^{2}$ be a connected non-degenerate submanifold in the $3$-dimensional Lorentz-Minkowski space $E_{1}^{3}$ and $\boldsymbol{x}:M^{2}\rightarrow E_{1}^{3}$ be a parametric representation of a surface in the Lorentz-Minkowski 3-space $E_{1}^{3}$ equipped with the induced metric. Let $(x, y, z)$ be a rectangular coordinate system of $E_{1}^{3}$. By saying Lorentz-Minkowski space $E_{1}^{3}$, we mean the Euclidean space $E^{3}$ with the standard metric given by
\begin{equation}
ds^{2} = -dx^{2}+ dy^{2}+dz^{2}.  \notag
\end{equation}

As well known that a submanifold is called a \emph{$k$-type} submanifold if its position vector $\mathbf{x}$ can be written as a sum of eigenvectors of the Laplace-Beltrami operator, $\Delta$, according to $k$ distinct eigenvalues, i.e.,
$\mathbf{x}=\mathbf{y}_0+\mathbf{y}_1+\dots+\mathbf{y}_k,$ for a constant vector $\mathbf{y}_0$ and smooth non-constant functions $\mathbf{y}_k$, $(i=1,\ldots,k)$ such that $\Delta \mathbf{y}_i =\lambda_i \mathbf{y}_i$, $\lambda_i \in \mathbb{R}$. The year 1966 was the beginning when Takahashi in \cite{T1} stated that spheres and minimal surfaces are the only ones in $E^{3}$ whose position vector $\boldsymbol{x}$ satisfies the relation
\begin{equation} \label{1.1}
\Delta^{I}\boldsymbol{x} =\lambda\boldsymbol{x}, \ \ \ \lambda \in R,
\end{equation}
where $\Delta^{I}$ is the Laplace operator associated with the $1^{st}$ fundamental form $I$ of the surface. Since the coordinate functions of $\boldsymbol{x}$ can be denoted as $(x_{1}, x_{2}, x_{3})$, then Takahashi's condition (\ref{1.1}) becomes

Let $(x_{1}, x_{2}, x_{3})$ be the component functions of $\boldsymbol{x}$. Then it is well-known that
\begin{equation} \label{1.2}
\Delta^{I}\boldsymbol{x} = (\Delta^{I}x_{1}, \Delta^{I}x_{2}, \Delta^{I}x_{3}).
\end{equation}

Thus Takahashi's condition (\ref{1.1}) becomes
\begin{equation} \label{1.3}
\Delta^{I}x_{i} =\lambda x_{i}, \ \ \ \ i = 1, 2, 3.
\end{equation}

Later, in \cite{G2} O. Garay generalized Takahashi's condition (\ref{1.3}). Actually, he studied surfaces of revolution in $E^{3}$, whose
component functions satisfy the condition
\begin{equation*} \label{1.4}
\Delta^{I}x_{i} =\lambda_{i} x_{i}, \ \ \ \ i = 1, 2, 3,
\end{equation*}
that is, the component functions are eigenfunctions of their Laplacian but not necessary with the same eigenvalue. Another generalization is to study surfaces whose position vector $\boldsymbol{x}$ satisfies a relation of the form
\begin{equation} \label{1.5}
\Delta ^{I}\boldsymbol{x} =A\boldsymbol{x},
\end{equation}
where $A\in \mathbb{Re}^{3\times 3}$.

Many results concerning this can be found in (\cite{B1}, \cite{B5}, \cite{B6}, \cite{D1}, \cite{D2}, \cite{G2}).
This type of study can be also extended to any smooth map, not necessary for the position vector of the surface, for example, the Gauss map of a surface. Regarding this see (\cite{A8}, \cite{A13}, \cite{A15}, \cite{G2}, \cite{B3}, \cite{B4}, \cite{B5}, \cite{B6}, \cite{D4}, \cite{G3}). Similarly, another extension can be drawn by applying the conditions stated before but for the second or third fundamental form of a surface \cite{S1}. Here again, many results can be found in (\cite{A1}, \cite{A2}, \cite{A14}, \cite{K8}, \cite{S0}, \cite{S2}).

On the other hand, all the ideas mentioned above can be applied in the Lorentz-Minkowski space $E_{1}^{3}$. So, an interesting geometric question has posed classify all the surfaces in $E_{1}^{3}$, which satisfy the condition
\begin{equation} \label{1.6}
\Delta ^{J}\boldsymbol{x} =A\boldsymbol{x},  \ \ \ J =I, II, III,
\end{equation}%
where $A\in \mathbb{Re}^{3\times 3}$ and $\Delta ^{J}$ is the Laplace operator, with respect to the fundamental form $J$.

Kaimakamis and Papantoniou in \cite{K0} solved the above question for the class of surfaces of revolution with respect to the second fundamental form, while in \cite{B7} Bekkar and Zoubir studied the same class of surfaces with respect to the first fundamental form satisfying
\begin{equation*}
\Delta x^{i} =\lambda^{i}x^{i},  \ \ \ \lambda^{i} \in R.
\end{equation*}


\section{Basic concepts}
Let $C: \boldsymbol{r}(s): s\in(a, b)\subset E \longrightarrow  E^{2}$ be a curve in a plane $E^{2}$ of $E_{1}^{3}$ and $l$ be a straight line of $E^{2}$ which does not intersect the curve $C$. A surface of revolution $M^{2}$ in $E_{1}^{3}$ is defined to be a non-degenerate surface, revolving the curve $C$ around the axis $l$. If the axis $l$ is timelike, then we consider that the $z$-axis as axis of revolution. If the axis $l$ is spacelike, then we may assume that the $x$-axis or $y$-axis as axis of revolution. Without loss of generality, we may consider the $x$-axis as the axis of revolution. If the axis is null, then we may assume that this axis is the line spanned by the vector (0,1,1) of the $yz$-plane.

We consider the axis of revolution is the $x$-axis (spacelike) and the curve $C$ is lying in the $xy$-plane. Then a parameterization of $C$ with respect to its arclength is $\boldsymbol{r}(s) = (f(s), g(s), 0)$ where $f, g$ are smooth functions. Without loss of generality, we may assume that $f(s) > 0, s \in(a, b)$. A surface of revolution $M^{2}$ in $E_{1}^{3}$ in a system of local curvilinear coordinates $(s, \theta)$ is given by:
\begin{equation}  \label{2.1}
\boldsymbol{x}(s, \theta) =\big(f(s)\cosh \theta, g(s),f(s)\sinh \theta \big).
\end{equation}
or
\begin{equation*}
\boldsymbol{x}(s, \theta) =\big(f(s)\sinh\theta, g(s),f(s)\cosh\theta\big).
\end{equation*}
In the case that the axis of revolution is the $z$-axis (timelike) and the curve $C$ is given $\boldsymbol{r}(s) = (f(s), 0, g(s))$ and lies in the $xz$-plane, the surface of revolution $M^{2}$ is given by:
\begin{equation}  \label{2.2}
\boldsymbol{x}(s, \theta) =\big(f(s)\cos\theta, f(s)\sin\theta, g(s)\big).
\end{equation}

Finally, if the axis of revolution is the line spanned by the vector (0, 1, 1) and the curve $C$ lies in the $yz$-plane, then the surface of revolution $M^{2}$ can be parametrized as
\begin{equation}  \label{2.3}
\boldsymbol{x}(s, \theta) =\big(\theta h(s),g(s)+\frac{1}{2}\theta^{2}h(s),f(s)+\frac{1}{2}\theta^{2}h(s)\big),
\end{equation}
where $h(s)=f(s)-g(s)\neq 0$.

We denote by $g = g_{km}, b = b_{km}$ and $e = e_{km}, k,m =1,2$ the first, second and third fundamental forms of $M^{2}$ respectively, where we put

\begin{equation*}
g_{11}= E= <\boldsymbol{x_{s}}, \boldsymbol{x_{s}}>,\ \ \ g_{12}= F= <\boldsymbol{x_{s}}, \boldsymbol{x_{\theta}}>,\ \ \ g_{22}= G= <\boldsymbol{x_{\theta}}, \boldsymbol{x_{\theta}}>,
\end{equation*}

\begin{equation*}
b_{11}=  L = <\boldsymbol{x_{ss}}, \boldsymbol{n}>,\ \ \ b_{12} = M = <\boldsymbol{x_{s\theta}}, \boldsymbol{n}>,\ \ \ b_{22}= N= <\boldsymbol{x_{\theta\theta}}, \boldsymbol{n}>,
\end{equation*}

\begin{equation*}
e_{11} = \frac{EM^{2}-2FLM+GL^{2}}{EG-F^{2}}= <\boldsymbol{n_{s}}, \boldsymbol{n_{s}}>,
\end{equation*}
\begin{equation*}
e_{12} = \frac{EMN-FLN+GLM-FM^{2}}{EG-F^{2}}= <\boldsymbol{n_{s}}, \boldsymbol{n_{\theta}}>,
\end{equation*}
\begin{equation*}
e_{22} = \frac{GM^{2}-2FNM+EN^{2}}{EG-F^{2}}= <\boldsymbol{n_{\theta}}, \boldsymbol{n_{\theta}}>,
\end{equation*}
which are the coefficients of the first, second, third fundamental form respectively, and $<,>$ is the Lorentzian metric.

For a sufficient differentiable function $p(u^{1}, u^{2})$ on $M^{2}$ The second Laplace operator according to the fundamental form $III$ of $M^{2}$ is defined by \cite{A4}.

\begin{equation*}  \label{2.5}
\Delta ^{III}p =-\frac{1}{\sqrt{e}}(\sqrt{e}e^{km}p_{/k})_{/m}.
\end{equation*}
where $p_{/k}:=\frac{\partial p}{\partial u^{k}}$, $e^{km}$ denote the components of the inverse tensor of $e_{km}$ and $e =\det(e_{km})$.
After a long computation, we arrive at

\begin{eqnarray}  \label{2.6}
\Delta ^{III}p&=& -\frac{\sqrt{\mid EG-F^{2}\mid}}{LN-M^{2}}\bigg(\bigg(\frac{(GM^{2}-2FNM+EN^{2})\frac{\partial p}{\partial s}}{(LN-M^{2})\sqrt{\mid EG-F^{2}\mid}}
  \notag \\
 && -\frac{(EMN-FLN+GLM-FM^{2})\frac{\partial p}{\partial\theta}}{(LN-M^{2})\sqrt{\mid EG-F^{2}\mid}}\bigg)_{s}   \notag \\
 && -\bigg(\frac{(EMN-FLN+GLM-FM^{2})\frac{\partial p}{\partial s}}{(LN-M^{2})\sqrt{\mid EG-F^{2}\mid}}   \notag \\
 &&-\frac{(EM^{2}-2FLM+GL^{2})\frac{\partial p}{\partial \theta}}{(LN-M^{2})\sqrt{\mid EG-F^{2}\mid}}\bigg)_{\theta}\bigg).
\end{eqnarray}
Here we have $LN-M^{2}\neq 0$, since the surface has no parabolic points.

\section{Proof of the main results}
In this paragraph we classify the surfaces of revolution $M^{2}$ satisfying the relation (\ref{1.6}). We distinguish the following three types according to whether these surfaces are determined.

Type I. The parametric representation of $M^{2}$ is given by by (\ref{2.1}). Then
\begin{equation}  \label{3.1}
f'^{2}(s)+g'^{2}(s) = 1,
\end{equation}
where $\prime := \frac{d}{ds}$, from which we obtain that
\begin{equation}  \label{3.2}
E = 1, \ \ \ F = 0, \ \ \ G = -f^{2}
\end{equation}
and
\begin{equation}  \label{3.3}
L = -f'g'' + g'f'', \ \ \ M = 0, \ \ \ N = fg'.
\end{equation}

Denoting by $\kappa$ the curvature of the curve $C$ and $r_{1}, r_{_{2}}$ the principal radii of curvature of $M^{2}$. We have
\begin{equation*}  \label{3.4}
r_{1} = \kappa, \ \ \ \ r_{2}= \frac{g'}{f}
\end{equation*}
and
\begin{align}  \label{3.5}
K = r_{1}r_{2} =\frac{\kappa g'}{f}= -\frac{f''}{f},\ \ \ \ 2H = r_{1}+r_{2} = \kappa+ \frac{g'}{f},
\end{align}
which are the Gauss and mean curvature of $M^{2}$ respectively. Since the relation (\ref{3.1}) holds, there exists a smooth function $\varphi = \varphi(s)$ such that
\begin{equation*}  \label{3.6}
f'= \cos\varphi, \ \ \ g' = \sin\varphi ,
\end{equation*}
where $\varphi= \varphi(s)$. Then $\kappa = \varphi'$ and relations (\ref{3.3}), (\ref{3.5}) become
\begin{equation}  \label{3.7}
L = -\varphi', \ \ \ M = 0, \ \ \ N = f\sin\varphi,
\end{equation}
\begin{align*}  \label{3.8}
K  =\frac{\varphi' \sin\varphi}{f}= -\frac{f''}{f},\ \ \ \ 2H = -\varphi'- \frac{\sin\varphi}{f}.
\end{align*}

We put $r = \frac{1}{r_{1}} + \frac{1}{r_{2}} = \frac{2H}{K}$. Thus we have
\begin{equation}  \label{3.9}
r= -\bigg(\frac{1}{\varphi'}+\frac{f}{\sin\varphi}\bigg).
\end{equation}

Taking the derivative of last equation, we get
\begin{equation}  \label{3.10}
r'= \bigg(\frac{\varphi''}{\varphi'^{2}}+\frac{f\varphi'\cos\varphi}{\sin^{2}\varphi}-\frac{\cos\varphi}{\sin\varphi}\bigg).
\end{equation}

From (\ref{2.6}), (\ref{3.2}) and (\ref{3.7}) we have
\begin{eqnarray}  \label{3.11}
\Delta ^{III}p= -\frac{1}{\varphi'^{2}}\frac{\partial^{2} p}{\partial s^{2}}+\frac{1}{\sin^{2}\varphi}\frac{\partial^{2}p}{\partial\theta^{2}}+\bigg(\frac{\varphi''}{\varphi'^{3}}-\frac{\cos\varphi}{\varphi'\sin\varphi}\bigg)\frac{\partial p}{\partial s}.
\end{eqnarray}

Let $(x_{1}, x_{2}, x_{3})$ be the coordinate functions of the position vector $\boldsymbol{x}$ of (\ref{2.3}). Then according to relation (\ref{1.2}), (\ref{3.11}) and taking into account (\ref{3.9}) and (\ref{3.10}) we find that
\begin{equation}  \label{3.12}
\Delta^{III}x_{1} =\Delta^{III}f(s)\cosh \theta = \bigg(-r\sin\varphi +r'\frac{\cos\varphi}{\varphi'}\bigg)\cosh \theta,
\end{equation}
\begin{equation}  \label{3.13}
\Delta^{III}x_{2} =\Delta^{III}g(s) = r\cos\varphi +r'\frac{\sin\varphi}{\varphi'},
\end{equation}
\begin{equation}  \label{3.14}
\Delta^{III}x_{3} =\Delta^{III}f(s)\sinh \theta = \bigg(-r\sin\varphi +r'\frac{\cos\varphi}{\varphi'}\bigg)\sinh \theta.
\end{equation}

We denote by $a_{ij},i,j=1,2,3,$ the entries of the matrix $A$, where all entries are real numbers. By using (\ref{3.12}), (\ref{3.13}) and (\ref{3.14}) condition (\ref{1.6}) is found to be equivalent to the following system
\begin{equation}  \label{3.15}
\bigg(-r\sin\varphi +r'\frac{\cos\varphi}{\varphi'}\bigg)\cosh \theta = a_{11}f(s)\cosh \theta +a_{12}g(s)+ a_{13}f(s)\sinh \theta,
\end{equation}
\begin{equation}  \label{3.16}
r\cos\varphi +r'\frac{\sin\varphi}{\varphi'} = a_{21}f(s)\cosh \theta +a_{22}g(s)+a_{23}f(s)\sinh \theta,
\end{equation}
\begin{equation}  \label{3.17}
\bigg(-r\sin\varphi +r'\frac{\cos\varphi}{\varphi'}\bigg)\sinh \theta = a_{31}f(s)\cosh \theta +a_{32}g(s)+a_{33}f(s)\sinh \theta.
\end{equation}

From (\ref{3.16}) it can be easily verified that $a_{21}= a_{23}= 0$. On the other hand, differentiating (\ref{3.15}) and (\ref{3.17}) twice with respect to $\theta$ we get that $a_{12}=a_{32}= 0$. So, the system is reduced to
\begin{equation}  \label{3.18}
\bigg(-r\sin\varphi +r'\frac{\cos\varphi}{\varphi'}\bigg)\cosh \theta = a_{11}f(s)\cosh \theta + a_{13}f(s)\sinh \theta,
\end{equation}
\begin{equation}  \label{3.19}
r\cos\varphi +r'\frac{\sin\varphi}{\varphi'} = a_{22}g(s),
\end{equation}
\begin{equation}  \label{3.20}
\bigg(-r\sin\varphi +r'\frac{\cos\varphi}{\varphi'}\bigg)\sinh \theta = a_{31}f(s)\cosh \theta +a_{33}f(s)\sinh \theta.
\end{equation}

But $\sinh \theta$ and $\cosh \theta$ are linearly independent functions of $\theta$, so we deduce that $a_{13}=a_{31}=0,a_{11}=a_{33}$. Putting $a_{11}=a_{33}=\lambda$ and $a_{22}=\mu$, we see that the system of equations (\ref{3.18}), (\ref{3.19}) and (\ref{3.20}) reduces now to the following two equations
\begin{equation}  \label{3.21}
-r\sin\varphi +r'\frac{\cos\varphi}{\varphi'} = \lambda f,
\end{equation}
\begin{equation}  \label{3.22}
r\cos\varphi +r'\frac{\sin\varphi}{\varphi'} = \mu g.
\end{equation}

Hence the matrix $A$ for which relation (\ref{1.6}) is satisfied becomes
\begin{equation*}
A =\left[
\begin{array}{ccc}
\lambda & 0 & 0 \\
0 & \mu & 0 \\
0 & 0& \lambda
\end{array}%
\right].
\end{equation*}

Solving the system (\ref{3.21}) and (\ref{3.22}) with respect to $r$ and $r'$, we conclude that
\begin{equation}  \label{3.23}
r'= \varphi'(\lambda f\cos\varphi + \mu g\sin\varphi) ,
\end{equation}
\begin{equation}  \label{3.24}
r = \mu g\cos\varphi -\lambda f\sin\varphi.
\end{equation}

Taking the derivative of (\ref{3.24}), we find
\begin{equation}  \label{3.25}
r'= \frac{1}{2}(\mu -\lambda)\cos\varphi \sin\varphi.
\end{equation}

We distinguish now the following cases:

\textit{Case I.} $\mu =\lambda = 0$. Thus, according to (\ref{3.24}) we have $r = 0$. Consequently $H =0$. Therefore $M^{2}$ is minimal and the corresponding matrix $A$ is the zero matrix.

\textit{Case II.} $\mu =\lambda \neq 0$. Then from (\ref{3.25}) we have $r' = 0$. If $\varphi'=0$, then $M^{2}$ would consist only of parabolic points, which has been excluded. Therefore we find that
\begin{equation*}  \label{3.26}
f(s)\cos\varphi + g(s)\sin\varphi= 0,
\end{equation*}
or
\begin{equation*}  \label{3.27}
ff' + gg'= 0.
\end{equation*}

Then $f^{2} + g^{2} = c^{2}, c\in R$ and $M^{2}$ is obviously satisfies the equation $x^{2}+y^{2}-z^{2}=c^{2}$ which is the pseudosphere $S^{2}_{1}(c)$ of $E_{1}^{3}$.

\textit{Case III.} $\lambda \neq 0, \mu =0$. Then the system (\ref{3.21}), (\ref{3.22}) is equivalently reduced to
\begin{equation*}  \label{3.28}
-r\sin\varphi +r'\frac{\cos\varphi}{\varphi'} = \lambda f(s),
\end{equation*}
\begin{equation*}  \label{3.29}
r\cos\varphi +r'\frac{\sin\varphi}{\varphi'} = 0.
\end{equation*}

From (\ref{3.24}) we have
\begin{equation}  \label{3.30}
r +\lambda f\sin\varphi= 0.
\end{equation}

On differentiating (\ref{3.30}) and taking into account (\ref{3.23}) with $\mu =0$, we obtain
\begin{equation*}  \label{3.31}
\lambda f\varphi'\cos\varphi +\lambda \cos\varphi\sin\varphi +\lambda f\varphi'\cos\varphi = 0
\end{equation*}
or
\begin{equation*}  \label{3.32}
\varphi' = -\frac{\sin\varphi}{2f}.
\end{equation*}

From (\ref{3.30}), (\ref{3.9}) and the last equation, we get
\begin{equation*}
\frac{f}{\sin\varphi}+\lambda f\sin\varphi= 0
\end{equation*}
or
\begin{equation*}  \label{3.33}
f(1+\lambda \sin^{2}\varphi)= 0.
\end{equation*}

A contradiction. Hence, there is no surface of revolution with parametric representation (\ref{2.1}) of $E_{1}^{3}$ satisfying (\ref{1.6}).

\textit{Case IV.} $\lambda=0, \mu \neq 0$. Then equations (\ref{3.21}), (\ref{3.22}) reduced to
\begin{equation*}  \label{3.34}
-r\sin\varphi +r'\frac{\cos\varphi}{\varphi'} = 0,
\end{equation*}
\begin{equation}  \label{3.35}
r\cos\varphi +r'\frac{\sin\varphi}{\varphi'} = \mu g.
\end{equation}

From (\ref{3.24}) we have
\begin{equation}  \label{3.36}
r - \mu g\cos\varphi = 0.
\end{equation}

Taking the derivative of (\ref{3.36}) and taking into account (\ref{3.23}) with $\lambda =0$, we find
\begin{equation*}  \label{3.37}
\mu g\varphi' \sin\varphi - \mu \cos\varphi\sin\varphi +\mu g\varphi'\sin \varphi = 0
\end{equation*}
or
\begin{equation}  \label{3.38}
\varphi' = \frac{\cos\varphi}{2g}.
\end{equation}

Taking the derivative of (\ref{3.38}), we find
\begin{equation}  \label{3.39}
3\varphi'\sin\varphi + 2g\varphi'' =0.
\end{equation}

On account of (\ref{3.35}), (\ref{3.9}) and (\ref{3.10}) it is easily verified that
\begin{equation}  \label{3.40}
\varphi'' =\frac{\varphi'^{2}}{\sin\varphi}(\mu g\varphi'+ 2\cos\varphi) .
\end{equation}

Inserting (\ref{3.38}) and (\ref{3.40}) in (\ref{3.39}) we conclude
\begin{equation*}  \label{3.41}
3 + (\frac{1}{2}\mu -1)\cos^{2}\varphi= 0 .
\end{equation*}

Here we have also a contradiction.

\textit{Case V.} $\lambda \neq 0, \mu \neq 0$.
We write equations (\ref{3.21}) and (\ref{3.22}) as follows
\begin{equation}  \label{3.42}
\frac{\sin\varphi}{\varphi'} + \frac{f}{\sin^{2}\varphi}+ \frac{\varphi''\cos\varphi}{\varphi'^{3}}-\frac{\cos^{2}\varphi}{\varphi'\sin\varphi} - \lambda f= 0,
\end{equation}
\begin{equation}  \label{3.43}
\frac{\varphi''\sin\varphi}{\varphi'^{3}} - \frac{2\cos\varphi}{\varphi'} - \mu g=0.
\end{equation}

From (\ref{3.43}) we have relation (\ref{3.40}). By eliminating $\varphi''$ from (\ref{3.42}) we get
\begin{equation}  \label{3.44}
\frac{1}{\varphi'\sin\varphi} + \frac{\mu g\cos\varphi}{\sin\varphi}+ \frac{f}{\sin^{2}\varphi} - \lambda f= 0.
\end{equation}

On differentiating the last equation and using (\ref{3.40}) we find 
\begin{equation}  \label{3.45}
\frac{2\mu g\varphi'}{\sin^{2}\varphi} + \frac{2f\varphi'\cos\varphi}{\sin^{3}\varphi} +\frac{2\cos\varphi}{\sin^{2}\varphi} -(\mu - \lambda)\cos\varphi =0.
\end{equation}

Multiplying (\ref{3.44}) by $\frac{2\varphi'}{\sin\varphi}$ and (\ref{3.45}) by $-\cos\varphi$ we obtain
\begin{equation}  \label{3.46}
\frac{2}{\sin^{2}\varphi} + \frac{2\mu g\varphi'\cos\varphi}{\sin^{2}\varphi}+ \frac{2f\varphi'}{\sin^{3}\varphi} - \frac{2\lambda\varphi' f}{\sin\varphi}= 0,
\end{equation}
\begin{equation}  \label{3.47}
-\frac{2\mu g\varphi'\cos\varphi}{\sin^{2}\varphi} - \frac{2f\varphi'\cos^{2}\varphi}{\sin^{3}\varphi} -\frac{2\cos^{2}\varphi}{\sin^{2}\varphi} +(\mu - \lambda)\cos^{2}\varphi =0.
\end{equation}

Combining (\ref{3.46}) and (\ref{3.47}) we conclude that
\begin{equation}  \label{3.48}
(\mu - \lambda)\cos^{2}\varphi -2(\lambda-1) \frac{f\varphi'}{\sin\varphi} +2=0
\end{equation}
or
\begin{equation*}  \label{3.49}
\frac{(\mu - \lambda)\cos^{2}\varphi}{\varphi'} -2(\lambda-1) \frac{f}{\sin\varphi} +\frac{2}{\varphi'}=0.
\end{equation*}
Taking the derivative of the above equation and using (\ref{3.25}) and (\ref{3.40}) we find
\begin{equation}  \label{3.50}
2(\mu+1)\cos\varphi + (\mu - \lambda)\mu g\varphi'\cos^{2}\varphi-2(\lambda-1) \frac{f\varphi'\cos\varphi} {\sin\varphi} +2\mu g\varphi'=0.
\end{equation}

Multiplying (\ref{3.48}) by $-\cos\varphi$, and adding the resulting equation to (\ref{3.50}) we get
\begin{equation*}  \label{3.51}
2\mu\cos\varphi +(2+ (\mu - \lambda)\cos^{2}\varphi)\mu g\varphi'-(\mu -\lambda)\cos^{3}\varphi =0
\end{equation*}
or
\begin{equation}  \label{3.52}
2\mu\cos^{2}\varphi +(2+ (\mu - \lambda)\cos^{2}\varphi)\mu g\varphi'\cos\varphi-(\mu -\lambda)\cos^{4}\varphi =0.
\end{equation}
On account of (\ref{3.44}) we find
\begin{equation}  \label{3.53}
\mu g\varphi'\cos\varphi= \lambda f\varphi'\sin\varphi - \frac{f \varphi'}{\sin\varphi} -1.
\end{equation}

Eliminating $\mu g\varphi'\cos\varphi$ from (\ref{3.52}) by using (\ref{3.53}), equation (\ref{3.52}) reduces to
\begin{eqnarray}  \label{3.54}
&& 2\mu\cos^{2}\varphi -(\mu -\lambda)\cos^{4}\varphi+  \notag \\
&&\big(2+(\mu - \lambda)\cos^{2}\varphi\big)\big((\lambda \sin^{2}\varphi -1) \frac{f \varphi'}{\sin\varphi} -1\big) =0.
\end{eqnarray}

But from (\ref{3.48}) we have
\begin{equation}  \label{3.55}
\frac{f\varphi'}{\sin\varphi}=\frac{(\mu - \lambda)\cos^{2}\varphi +2}{2(\lambda-1)}.
\end{equation}

Obviously $\lambda \neq 1$ because otherwise, from (\ref{3.48}) we would have
\begin{equation*}  \label{3.56}
(\mu - \lambda)\cos^{2}\varphi +2=0.
\end{equation*}

A contradiction. Now, by inserting (\ref{3.55}) in (\ref{3.54}) we obtain
\begin{eqnarray*}  \label{3.57}
&& -\lambda(\mu - \lambda)^{2}\cos^{4}\varphi +(\mu -\lambda)\big((\mu - \lambda)(\lambda-1) - 6\lambda +2\big)\cos^{2}\varphi \notag \\
&&+6\mu(\lambda-1)- 2\lambda (\lambda+1) =0.
\end{eqnarray*}

This relation, however, is valid for a finite number of values of $\varphi$. So in this case there are no surfaces of revolution with the required property. So we proved the following
\begin{theorem}\label{T9}
Let $\boldsymbol{x}:M^{2} \longrightarrow E^{3}_{1}$ be a surface of revolution given by (\ref{2.1}). Then $\boldsymbol{x}$ satisfies (\ref{1.6}) regarding to the third fundamental form if and only if the following statements hold true
\begin{itemize}
\item $M^{2}$ is the pseudosphere $S^{2}_{1}(c)$ of real radius $c$,
\item $M^{2}$ has zero mean curvature.
\end{itemize}
\end{theorem}

Type II. The parametric representation of $M^{2}$ is given by (\ref{2.2}). Then the tangent vector of the revolving curve is
\begin{equation*}  \label{3.58}
<\boldsymbol{x}', \boldsymbol{x}'> = f'^{2}-g'^{2}= \pm 1.
\end{equation*}

We assume that
\begin{equation}  \label{3.59}
 f'^{2}-g'^{2}= 1, \ \ \ \forall s \in(a, b).
\end{equation}

Then the components of the first and second fundamental forms are respectively
\begin{equation*}  \label{3.60}
E = 1, \ \ \ F = 0, \ \ \ G =f^{2},
\end{equation*}

\begin{equation*}  \label{3.61}
L =f'g'' - g'f'', \ \ \ M = 0, \ \ \ N = fg'.
\end{equation*}

From the equation (\ref{3.59}) it is obviously clear that there exist a smooth function $\varphi = \varphi(s)$ such that
\begin{equation*}  \label{3.62}
f'= \cosh\varphi, \ \ \ g' = \sinh\varphi.
\end{equation*}

On the other hand we have
\begin{equation*}  \label{3.62}
r_{1} = \kappa = \varphi', \ \ \ \ r_{2}= \frac{g'}{f} = \frac{\sinh\varphi}{f}
\end{equation*}
and
\begin{align*}  \label{3.63}
K = r_{1}r_{2} =\frac{\kappa g'}{f}= -\frac{f''}{f}= \frac{\varphi'\sinh\varphi}{f},\ \ \ \ 2H = r_{1}+r_{2} = \varphi'+ \frac{\sinh\varphi}{f}.
\end{align*}
Here we have
\begin{equation}  \label{3.64}
r= \frac{1}{\varphi'}+\frac{f}{\sinh\varphi}.
\end{equation}

Taking the derivative of last equation, we get
\begin{equation}  \label{3.65}
r'= -\frac{\varphi''}{\varphi'^{2}}-\frac{f\varphi'\cosh\varphi}{\sinh^{2}\varphi}+\frac{\cosh\varphi}{\sinh\varphi}.
\end{equation}

On the other hand
\begin{eqnarray}  \label{3.66}
\Delta ^{III}p= -\frac{1}{\varphi'^{2}}\frac{\partial^{2} p}{\partial s^{2}}-\frac{1}{\sinh^{2}\varphi}\frac{\partial^{2}p}{\partial\theta^{2}}+\bigg(\frac{\varphi''}{\varphi'^{3}}-\frac{\cosh\varphi}{\varphi'\sinh\varphi}\bigg)\frac{\partial p}{\partial s}.
\end{eqnarray}

According to relation (\ref{2.2}) and (\ref{3.66}) we find that
\begin{equation*}  \label{3.67}
\Delta^{III}x_{1} =\Delta^{III}f(s)\cos\theta = \bigg(-r\sinh\varphi -r'\frac{\cosh\varphi}{\varphi'}\bigg)\cos\theta,
\end{equation*}
\begin{equation*}  \label{3.68}
\Delta^{III}x_{2} =\Delta^{III}f(s)\sin\theta = \bigg(-r\sinh\varphi -r'\frac{\cosh\varphi}{\varphi'}\bigg)\sin\theta,
\end{equation*}
\begin{equation*}  \label{3.69}
\Delta^{III}x_{3} =\Delta^{III}g(s) = -r\cosh\varphi -r'\frac{\sinh\varphi}{\varphi'}.
\end{equation*}

Let now $\Delta ^{III} \boldsymbol{x} = A\boldsymbol{x}$. Thus, as in the former paragraph, we find
\begin{equation*}  \label{3.70}
\bigg(-r\sinh\varphi -r'\frac{\cosh\varphi}{\varphi'}\bigg)\cos\theta = a_{11}f(s)\cos\theta +a_{12}f(s)\sin\theta+ a_{13}g(s),
\end{equation*}
\begin{equation*}  \label{3.71}
\bigg(-r\sinh\varphi -r'\frac{\cosh\varphi}{\varphi'}\bigg)\sin\theta = a_{21}f(s)\cos\theta +a_{22}f(s)\sin\theta + a_{23}g(s),
\end{equation*}
\begin{equation*}  \label{3.72}
-r\cosh\varphi -r'\frac{\sinh\varphi}{\varphi'}= a_{31}f(s)\cos\theta +a_{32}f(s)\sin\theta+a_{33}g(s).
\end{equation*}

Applying the same algebraic methods, used in the previous type, this system of equations reduced to
\begin{equation}  \label{3.73}
-r\sinh\varphi -r'\frac{\cosh\varphi}{\varphi'} = \lambda f,
\end{equation}
\begin{equation}  \label{3.74}
-r\cosh\varphi -r'\frac{\sinh\varphi}{\varphi'}= \mu g,
\end{equation}
where $a_{11} = a_{22} =  \lambda, a_{33} = \mu, \lambda,\mu \in R$.
Solving the system (\ref{3.73}) and (\ref{3.74}) with respect to $r$ and $r'$, we conclude that
\begin{equation}  \label{3.75}
r'= \varphi'(-\lambda f\cosh\varphi + \mu g\sinh\varphi) ,
\end{equation}
\begin{equation}  \label{3.76}
r = \lambda f\sinh\varphi -\mu g\cosh\varphi.
\end{equation}

Similarly, we have the following five cases according to the values of $\lambda,\mu$.

\textit{Case I.} $\lambda = \mu =0$. Thus from (\ref{3.76}) we conclude that $r = 0$. Consequently $H =0$. Therefore $M^{2}$ is minimal and the corresponding matrix $A$ is the zero matrix.

\textit{Case II.} $\mu =\lambda \neq 0$. Then from (\ref{3.75}) we have $r' = 0$. If $\varphi'=0$, then $M^{2}$ would consist only of parabolic points, which has been excluded. Therefore we find that
\begin{equation*}  \label{3.77}
-f\cosh\varphi + g\sinh\varphi= 0,
\end{equation*}
or
\begin{equation*}  \label{3.78}
-ff' + gg'= 0.
\end{equation*}

Then $g^{2} -f^{2}= \pm c^{2}, c\in R$

and therefore, $M^{2}$ is obviously either the pseudosphere $S^{2}_{1}(c)$ of real radius $c$, given by the equation $x^{2}+y^{2}-z^{2}=c^{2}$, or the pseudosphere $H^{2}(c)$ with imaginary radius, given by $x^{2}+y^{2}-z^{2}= -c^{2}$.

\textit{Case III.} $\lambda \neq 0, \mu =0$. Then the system (\ref{3.73}), (\ref{3.74}) is reduced to
\begin{equation*}  \label{3.79}
-r\sinh\varphi -r'\frac{\cosh\varphi}{\varphi'} = \lambda f(s),
\end{equation*}
\begin{equation*}  \label{3.80}
r\cosh\varphi +r'\frac{\sinh\varphi}{\varphi'} = 0.
\end{equation*}

From (\ref{3.76}) we have
\begin{equation}  \label{3.81}
r -\lambda f\sinh\varphi = 0.
\end{equation}

On differentiating (\ref{3.81}) and taking into account (\ref{3.75}) with $\mu =0$, we obtain
\begin{equation*}  \label{3.82}
\varphi' = -\frac{\sinh\varphi}{2f}.
\end{equation*}

From (\ref{3.64}), (\ref{3.76}) and the last equation, we get

\begin{equation*}  \label{3.83}
f(1-\lambda \sin^{2}\varphi)= 0.
\end{equation*}

A contradiction. Hence, there are no surfaces of revolution with parametric representation (\ref{2.2}) of $E_{1}^{3}$ satisfying (\ref{1.6}).

\textit{Case IV.} $\lambda=0, \mu \neq 0$. Then equations (\ref{3.73}), (\ref{3.74}) reduced to
\begin{equation*}  \label{3.84}
-r\sinh\varphi -r'\frac{\cosh\varphi}{\varphi'} = 0,
\end{equation*}
\begin{equation}  \label{3.85}
-r\cosh\varphi -r'\frac{\sinh\varphi}{\varphi'} = \mu g.
\end{equation}

From (\ref{3.76}) we have
\begin{equation}  \label{3.86}
r + \mu g\cosh\varphi = 0.
\end{equation}

Taking the derivative of (\ref{3.86}) and taking into account (\ref{3.75}) with $\lambda =0$, we find
\begin{equation}  \label{3.87}
\varphi' = -\frac{\cosh\varphi}{2g}.
\end{equation}

Taking the derivative of (\ref{3.87}), we find
\begin{equation}  \label{3.88}
3\varphi'\sinh\varphi + 2g\varphi'' =0.
\end{equation}

On account of (\ref{3.85}), (\ref{3.64}) and (\ref{3.65}) it is easily verified that
\begin{equation}  \label{3.89}
\varphi'' =\frac{\varphi'^{2}}{\sinh\varphi}(\mu g\varphi'+ 2\cosh\varphi) .
\end{equation}

Inserting (\ref{3.87}) and (\ref{3.89}) in (\ref{3.88}) we conclude
\begin{equation*}  \label{3.90}
3 - (\frac{1}{2}\mu +1)\cosh^{2}\varphi= 0,
\end{equation*}
which we easily get a contradiction.

\textit{Case V.} $\lambda \neq 0, \mu \neq 0$.
We write equations (\ref{3.73}) and (\ref{3.74}) as follows
\begin{equation}  \label{3.91}
-\frac{\sinh\varphi}{\varphi'} + \frac{f}{\sinh^{2}\varphi}+ \frac{\varphi''\cosh\varphi}{\varphi'^{3}}-\frac{\cosh^{2}\varphi}{\varphi'\sinh\varphi} - \lambda f= 0,
\end{equation}
\begin{equation}  \label{3.92}
\frac{\varphi''\sinh\varphi}{\varphi'^{3}} - \frac{2\cosh\varphi}{\varphi'} - \mu g=0.
\end{equation}

From (\ref{3.92}) we have relation (\ref{3.89}). By eliminating $\varphi''$ from (\ref{3.91}) we get
\begin{equation}  \label{3.93}
\frac{1}{\varphi'\sinh\varphi} + \frac{\mu g\cosh\varphi}{\sinh\varphi}+ \frac{f}{\sinh^{2}\varphi} - \lambda f= 0.
\end{equation}

On differentiating the last equation and using (\ref{3.89}) we find
\begin{equation}  \label{3.94}
\frac{2\mu g\varphi'}{\sinh^{2}\varphi} + \frac{2f\varphi'\cosh\varphi}{\sinh^{3}\varphi} +\frac{2\cosh\varphi}{\sinh^{2}\varphi} -(\mu - \lambda)\cosh\varphi =0.
\end{equation}

Multiplying (\ref{3.93}) by $\frac{2\varphi'}{\sinh\varphi}$ and (\ref{3.94}) by $-\cosh\varphi$ we obtain
\begin{equation}  \label{3.95}
\frac{2}{\sinh^{2}\varphi} + \frac{2\mu g\varphi'\cosh\varphi}{\sinh^{2}\varphi}+ \frac{2f\varphi'}{\sinh^{3}\varphi} - \frac{2\lambda\varphi' f}{\sinh\varphi}= 0,
\end{equation}
\begin{equation}  \label{3.96}
-\frac{2\mu g\varphi'\cosh\varphi}{\sinh^{2}\varphi} - \frac{2f\varphi'\cosh^{2}\varphi}{\sinh^{3}\varphi} -\frac{2\cosh^{2}\varphi}{\sinh^{2}\varphi} +(\mu - \lambda)\cosh^{2}\varphi =0.
\end{equation}

Combining (\ref{3.95}) and (\ref{3.96}) we conclude that
\begin{equation}  \label{3.97}
(\mu - \lambda)\cosh^{2}\varphi -2(\lambda+1) \frac{f\varphi'}{\sinh\varphi} -2=0
\end{equation}
or
\begin{equation*}  \label{3.98}
\frac{(\mu - \lambda)\cosh^{2}\varphi}{\varphi'} -2(\lambda+1) \frac{f}{\sinh\varphi} -\frac{2}{\varphi'}=0.
\end{equation*}
Taking the derivative of the above equation and using (\ref{3.89}) we find
\begin{equation}  \label{3.99}
2(\mu+1)\cosh\varphi + (\mu - \lambda)\mu g\varphi'\cosh^{2}\varphi-2(\lambda+1) \frac{f\varphi'\cosh\varphi} {\sinh\varphi} -2\mu g\varphi'=0.
\end{equation}

Multiplying (\ref{3.97}) by $-\cosh\varphi$, and adding the resulting equation to (\ref{3.99}) we get
\begin{equation*}  \label{3.100}
2(\mu+2)\cosh\varphi -(2- (\mu - \lambda)\cos^{2}\varphi)\mu g\varphi'-(\mu -\lambda)\cosh^{3}\varphi =0
\end{equation*}
or
\begin{equation}  \label{3.101}
2(\mu+2)\cosh^{2}\varphi -(2- (\mu - \lambda)\cosh^{2}\varphi)\mu g\varphi'\cos\varphi-(\mu -\lambda)\cos^{4}\varphi =0.
\end{equation}
On account of (\ref{3.93}) we find
\begin{equation}  \label{3.102}
\mu g\varphi'\cosh\varphi= \lambda f\varphi'\sinh\varphi - \frac{f \varphi'}{\sinh\varphi} -1.
\end{equation}

Eliminating $\mu g\varphi'\cosh\varphi$ from (\ref{3.101}) by using (\ref{3.102}), we get
\begin{eqnarray}  \label{3.103}
&& 2(\mu+2)\cosh^{2}\varphi -(\mu -\lambda)\cosh^{4}\varphi-  \notag \\
&&\big(2-(\mu - \lambda)\cosh^{2}\varphi\big)\big((\lambda \sinh^{2}\varphi -1) \frac{f \varphi'}{\sinh\varphi} -1\big) =0.
\end{eqnarray}

But from (\ref{3.97}) we have
\begin{equation}  \label{3.104}
\frac{f\varphi'}{\sinh\varphi}=\frac{2-(\mu - \lambda)\cosh^{2}\varphi}{2(\lambda+1)}.
\end{equation}

Obviously $\lambda \neq -1$ because otherwise, from (\ref{3.97}) we would have
\begin{equation*}  \label{3.105}
(\mu - \lambda)\cosh^{2}\varphi -2=0.
\end{equation*}

A contradiction. Now, by inserting (\ref{3.104}) in (\ref{3.103}) we obtain
\begin{eqnarray*}  \label{3.57}
&& -\lambda(\mu - \lambda)^{2}\cos^{6}\varphi +(\mu -\lambda)\big((\mu - \lambda)(\lambda-1) +4\lambda\big)\cos^{4}\varphi \notag \\
&&+(6\lambda^{2}-2\lambda -2\mu -2\lambda \mu +8)\cos^{2}\varphi+ 8(\lambda+1) =0.
\end{eqnarray*}

This relation, however, is valid for a finite number of values of $\varphi$. So in this case there are no surfaces of revolution with the required property. So we proved the following
\begin{theorem}\label{T10}
Let $\boldsymbol{x}:M^{2} \longrightarrow E^{3}_{1}$ be a surface of revolution given by (\ref{2.2}). Then $\boldsymbol{x}$ satisfies (\ref{1.6}) regarding to the third fundamental form if and only if the following statements hold true
\begin{itemize}
\item $M^{2}$ is the pseudosphere $S^{2}_{1}(c)$ of real or imaginary radius $c$,
\item $M^{2}$ has zero mean curvature.
\end{itemize}
\end{theorem}

Type III. The parametric representation of $M^{2}$ is given by (\ref{2.3}), i.e.,
\begin{equation*}
\boldsymbol{x}(s, \theta) =\big(f(s)+\frac{1}{2}\theta^{2}h(s),g(s)+\frac{1}{2}\theta^{2}h(s),\theta h(s)\big),
\end{equation*}
where $h(s)=f(s)-g(s)\neq 0$. Since, $M^{2}$ is non-degenerate, $f'(s)^2-g'(s)^2$ never vanishes, and so $h'(s)=f'(s)-g'(s)\neq 0$ everywhere. Now, we may take the parameter in such a way that
\begin{equation*}
    h(s)=-2s.
\end{equation*}

Assume that $k(s)=g(s)-s$, then
\begin{equation*}
    f(s)=k(s)-s \quad g(s)=k(s)+s,
\end{equation*}
(see for example, \cite{K1}). Therefore $M^{2}$ can be reparametrized as follows
\begin{equation}  \label{4.43}
\boldsymbol{x}(s, \theta) =\big(k -s - \theta^{2}s ,k +s-\theta^{2}s,-2s\theta\big),
\end{equation}
with the profile curve given in \eqref{pccase3} becomes
\begin{equation} \label{npccase3}
    \boldsymbol{r}(s)=(0,k(s)-s,k(s)+s).
\end{equation}

By using the tangent vector fields, $\boldsymbol{x}_s$ and $\boldsymbol{x}_{\theta}$ of $M^{2}$, we get the components of the first fundamental form of it as
\begin{equation*}  \label{4.3}
E = 4k'(s), \quad F = 0, \quad G =4s^{2}.
\end{equation*}

Now, let $M^{2}$ be spacelike surface, i.e., $k'(s)>0$. Then, the timelike unit normal vector field $\boldsymbol{N}$ of $M^{2}$ is given by
\begin{equation} \label{case3N}
    \boldsymbol{N}=\frac{1}{2\sqrt{k'}}(\theta^2+1,\theta^2-1,2\theta)+\frac{\sqrt{k'}}{2}(1,1,0).
\end{equation}

Then the components of the second fundamental forms are given by
\begin{equation*}  \label{4.4}
L = -\frac{k''}{\sqrt{k'}}, \ \ \ M = 0, \ \ \ N = \frac{2s}{\sqrt{k'}}.
\end{equation*}

Thus the relation (\ref{2.6}) becomes
\begin{equation}  \label{4.44}
\Delta ^{III}p= -\frac{4k'^{2}}{k''^{2}}\frac{\partial^{2} p}{\partial s^{2}}- k'\frac{\partial^{2}p}{\partial\theta^{2}}+\frac{2k'}{k''^{3}}\big(2k'k'''-k''^{2}\big)\frac{\partial p}{\partial s}.
\end{equation}

According to relations (\ref{2.3}) and (\ref{4.44}) we find that
\begin{equation*}  \label{4.10}
\Delta^{III}x_{1} =\Delta^{III}(k-s-s\theta^{2}) = \frac{2k'}{{k''}^3}\big(2k'k'''-{k''}^2\big)(k'-1-\theta^2)-\frac{4{k'}^2}{k''}+2sk',
\end{equation*}
\begin{equation*}  \label{4.9}
\Delta^{III}x_{2} =\Delta^{III}(k+s-s\theta^{2}) = \frac{2k'}{{k''}^3}\big(2k'k'''-{k''}^2\big)(k'+1-\theta^2)-\frac{4{k'}^2}{k''}+2sk',
\end{equation*}
\begin{equation*}  \label{4.8}
\Delta^{III}x_{3} =\Delta^{III}(-2s\theta) = -\frac{4k'}{{k''}^3}\big(2k'k'''-{k''}^2\big)\theta.
\end{equation*}

Let now $\Delta ^{III} \boldsymbol{x} = A\boldsymbol{x}$. Then
\begin{eqnarray}  \label{4.13}
&&\frac{2k'}{{k''}^3}\big(2k'k'''-{k''}^2\big)(k'-1-\theta^2)-\frac{4{k'}^2}{k''}+2sk'=   \notag \\
&& a_{11}(k-s-s\theta^{2}) +a_{12}(k+s-s\theta^{2})+a_{13}(-2s\theta),
\end{eqnarray}
\begin{eqnarray}  \label{4.12}
&&\frac{2k'}{{k''}^3}\big(2k'k'''-{k''}^2\big)(k'+1-\theta^2)-\frac{4{k'}^2}{k''}+2sk' =   \notag \\
&& a_{21}(k-s-s\theta^{2}) +a_{22}(k+s-s\theta^{2}) + a_{23}(-2s\theta),
\end{eqnarray}
\begin{equation}  \label{4.11}
-\frac{4k'}{{k''}^3}\big(2k'k'''-{k''}^2\big)\theta = a_{31}(k-s-s\theta^{2}) +a_{32}(k+s-s\theta^{2})+ a_{33}(-2s\theta).
\end{equation}

Regarding the above equations as polynomials in $\theta$, so from the coefficients of (\ref{4.11}) we get
\begin{equation}  \label{4.14}
(a_{31}+a_{32})s = 0,
\end{equation}
\begin{equation}  \label{4.15}
\frac{2k'}{{k''}^3}\big(2k'k'''-{k''}^2\big)=a_{33}s,
\end{equation}
\begin{equation}  \label{4.16}
(a_{32}-a_{31})s+(a_{31}+a_{32})k = 0.
\end{equation}

From the coefficients of (\ref{4.12}) we find
\begin{equation}  \label{4.17}
\frac{2k'}{{k''}^3}\big(2k'k'''-{k''}^2\big)=(a_{21}+a_{22})s,
\end{equation}
\begin{equation}  \label{4.18}
a_{23}s = 0,
\end{equation}
\begin{equation}  \label{4.19}
\frac{2k'}{{k''}^3}\big(2k'k'''-{k''}^2\big)(k'+1)-\frac{{4k'}^2}{k''}+2sk'=(a_{21}+a_{22})k+(a_{22}-a_{21})s.
\end{equation}

From the coefficients of (\ref{4.13}) we get
\begin{equation}  \label{4.20}
\frac{2k'}{{k''}^3}\big(2k'k'''-{k''}^2\big)=(a_{11}+a_{12})s,
\end{equation}
\begin{equation}  \label{4.21}
a_{13}s = 0,
\end{equation}
\begin{equation}  \label{4.22}
\frac{2k'}{{k''}^3}\big(2k'k'''-{k''}^2\big)(k'-1)-\frac{{4k'}^2}{k''}+2sk'=(a_{11}+a_{12})k+(a_{12}-a_{11})s.
\end{equation}

It is easily verified that
\begin{equation*}  \label{4.23}
a_{23} = a_{31} = a_{32} = a_{13} =0.
\end{equation*}

On the other hand, from (\ref{4.15}), (\ref{4.17}) and (\ref{4.20}) we find
\begin{equation*}  \label{4.24}
a_{11} + a_{12} = a_{33} = a_{21} + a_{22},
\end{equation*}
from which we obtain
\begin{equation}  \label{4.25}
a_{12} = a_{33} - a_{11}, \quad a_{21} = a_{33}- a_{22}.
\end{equation}

Moreover, by considering \eqref{4.15} and \eqref{4.25} in \eqref{4.19} and \eqref{4.22}, respectively, we get

\begin{equation} \label{4.190}
    a_{33}s(k'+1)-\frac{{4k'}^2}{k''}+2sk'=a_{33}k+(a_{22}-a_{21})s
\end{equation}
and
\begin{equation} \label{4.220}
    a_{33}s(k'-1)-\frac{{4k'}^2}{k''}+2sk'=a_{33}k+(a_{12}-a_{11})s.
\end{equation}

By subtracting (\ref{4.190}) from (\ref{4.220}) we obtain
\begin{equation}  \label{4.26}
(a_{11}- a_{21})s+(a_{22}- a_{12})s-2a_{33}s=0.
\end{equation}

From (\ref{4.25}) and (\ref{4.26}) we find
\begin{equation}  \label{4.27}
a_{21}=- a_{12}.
\end{equation}

Taking into account relations (\ref{4.26}) and (\ref{4.27}), we get
\begin{equation*}  \label{4.31}
a_{11}+a_{22} =2a_{33}.
\end{equation*}

We put $a_{11} = \lambda$ and $a_{22} = \mu$, so the matrix $A$ for which relation (\ref{1.6}) is satisfied takes finally the following form
\begin{equation*}
A =\left[
\begin{array}{ccc}
\lambda &\frac{1}{2}(\mu-\lambda) & 0 \\
\frac{1}{2}(\lambda-\mu)&\mu & 0\\
0 &0 & \frac{1}{2}(\lambda+\mu)
\end{array}%
\right].
\end{equation*}

Hence system of equations [(\ref{4.14}),... (\ref{4.22})] reduces to the following two equations
\begin{equation}  \label{4.45}
\frac{2k'}{k''^{3}}(2k'k'''-k''^{2})= a_{33}s,
\end{equation}
\begin{equation}  \label{4.46}
(a_{33}+2)k's + 2a_{12}s-\frac{4k'^{2}}{k''} - a_{33} k=0,
\end{equation}
where, as we mention before, $a_{33} =\frac{1}{2}(\lambda+\mu)$ and $a_{12} =\frac{1}{2}(\mu-\lambda)$.

Solving the system of equations (\ref{4.45}) and (\ref{4.46}) with respect to $\lambda$ and $\mu$ we find
\begin{equation}  \label{4.47}
\lambda = \frac{k'(2s-k+sk')}{s^{2}k''}\Big(\frac{2k'k'''}{{k''}^2}- 1\Big)- \frac{2{k'}^2}{sk''} +k',
\end{equation}
\begin{equation}  \label{4.48}
\mu = \frac{k'(2s+k-sk')}{s^{2}k''}\Big(\frac{2k'k'''}{{k''}^2}- 1\Big)+ \frac{2{k'}^2}{sk''} -k'.
\end{equation}

\textit{Case I.} $\lambda = \mu =0$. Thus from (\ref{4.47}) and \eqref{4.48} we conclude that $k = as^3+b$ with $a>0$, $b$ is a constant and $s\neq 0$. Consequently, $H =0$. Therefore $M^{2}$ is minimal and the corresponding matrix $A$ is the zero matrix.

\textit{Case II.} $\lambda=\mu \neq 0$. Thus from \textit{Case I}, $k \neq as^3+b$. Now, from (\ref{3.88}) we get $a_{23}=0$, and so
\begin{equation}
    \frac{(k-sk')(2k'k'''-{k''}^2)}{s^{2}{k''}^3}+\frac{2k'}{sk''} -1=0,
\end{equation}
whose solution is $\displaystyle k(s)=\pm \frac{c^2}{4s}$. By considering \eqref{npccase3}, we conclude $\boldsymbol{r}$ is a spherical curve and so the surface $M^{2}$ is an open piece of the pseudo-sphere $\mathbb{S}^2_1(0,c)$ or the hyperbolic space $\mathbb{H}^2(0,c)$.

\textit{Case III.} $\lambda \neq 0, \mu =0$. By considering the last assumption in \eqref{4.48}, i.e. $\mu=0$, we have
\begin{equation*}
   \frac{2k'}{sk''}\Big(\frac{2k'k'''}{{k''}^2}- 1\Big)=\frac{k'(-k+sk')}{s^{2}k''}\Big(\frac{2k'k'''}{{k''}^2}- 1\Big)- \frac{2{k'}^2}{sk''} +k'.
\end{equation*}
By substituting this into \eqref{4.47}, we get
\begin{equation*}
    \lambda=\frac{4k'}{sk''}\Big(\frac{2k'k'''}{{k''}^2}- 1\Big),
\end{equation*}
where $\lambda$ is non-zero function. Since there is no $k$ function to implement in both conditions, so there is no  surface of revolution that fulfills these conditions.

\textit{Case IV.} $\lambda=0, \mu \neq 0$. Similarly, we get a contradiction as in \textit{Case III.}

\textit{Case V.} $\lambda \neq \mu $ and $\lambda \neq 0, \mu \neq 0$. In this case, the above two relations \eqref{4.47} and \eqref{4.48} are valid only when $\lambda$ and $\mu$ are functions of $s$. Thus there are no surfaces of revolution with the required property. So we proved the following:
\begin{theorem}\label{T11}
Let $\boldsymbol{x}:M^{2} \longrightarrow E^{3}_{1}$ be a surface of revolution given by (\ref{2.3}). Then $\boldsymbol{x}$ satisfies (\ref{1.6}) regarding to the third fundamental form if and only if the following statements hold true
\begin{itemize}
\item $M^{2}$ has zero mean curvature,
\item $M^{2}$ is an open piece of the pseudo sphere $\mathbb{S}^{2}_{1}(0,c)$ of real radius $c$,
\item $M^{2}$ is an open piece of the hyperbolic space $\mathbb{H}^{2}_{1}(0,c)$ of real radius $c$.
\end{itemize}
\end{theorem}

Finally, as we know that the minimal surfaces of revolution with non-lightlike axis are congruent to a part of the catenoid and also with lightlike axis are congruent to a part of the surface of Enneper, (see for more details \cite{W1990}). Now, by combining Theorem 1, Theorem 2, Theorem 3 and \cite{W1990}:

\begin{theorem} (Classification)
Let $\boldsymbol{x}:M^{2} \longrightarrow E^{3}_{1}$ be a surface of revolution satisfying (\ref{1.6}) regarding the third fundamental form. Then $M$ is one of the following:
\begin{itemize}
    \item $M^{2}$ is an open part of catenoid of the 1st kind, the 2nd kind, the 3rd kind, the 4th kind or the 5th kind.
    \item $M^{2}$ is an open part of the surface of Enneper of the 2nd kind or the 3rd kind,
    \item $M^{2}$ is an open part of the pseudo sphere $\mathbb{S}^{2}_{1}(0,c)$ centered at the origin with radius $c$,
\item $M^{2}$ is an open part of the hyperbolic space $\mathbb{H}^{2}_{1}(0,c)$ centered at the origin with radius $c$.
\end{itemize}
\end{theorem}

\end{document}